\documentclass[a4paper,12pt]{article}
\usepackage{amsmath}
\usepackage{amssymb}
\usepackage{tabularx}
\usepackage{enumerate}

\newtheorem{thm}{Theorem}[section]
\newtheorem{lem}[thm]{Lemma}
\newtheorem{cor}[thm]{Corollary}

\newtheorem{rmk}{Remark}[section]
\newtheorem{defi}{Definition}[section]
\newtheorem{pppp}{Proof}

\newcommand{\qed}{\hspace{1em}\mbox{\raisebox{0.65ex}{\fbox{}}}}

\numberwithin{equation}{section}

\newcommand{\be}{\begin{equation}}
\newcommand{\ee}{\end{equation}}
\newcommand\bes{\begin{eqnarray}} \newcommand\ees{\end{eqnarray}}
\newcommand{\bess}{\begin{eqnarray*}}
\newcommand{\eess}{\end{eqnarray*}}

\newcommand{\R}{\mathbb{R}}
\newcommand{\bpf}{{\bf Proof:\ \ }}
\newcommand{\epf}{\mbox{}\hfill $\Box$}

\begin{document}

\thispagestyle{empty}

\title{The spreading fronts of an infective environment in a man-environment-man epidemic model\thanks{This work is supported by
Collaborative Research Project under the NRF-NSFC Cooperative Program (11211140236) and
the National Research Foundation of Korea (NRF) Grant (NRF-2012K1A2B1A03000598).
}}
\date{\empty}

\author{Inkyung Ahn$^{a}$ and Zhigui Lin$^b$\\
{\small $^a$Department of Mathematics, Korea University, Jochiwon,}\\
{\small Chung-Nam 339-700, South Korea}\\
{\small $^b$School of Mathematical Science, Yangzhou University, Yangzhou 225002, China}
}
 \maketitle

\begin{quote}
\noindent
{\bf Abstract.} { 
\small A reaction-diffusion model is investigated to understand infective environments in a man-environment-man epidemic model.
 The free boundary is introduced to
describe the expanding front of an infective environment induced by fecally-orally transmitted disease.
The basic reproduction number $R^F_0(t)$ for the free boundary problem is introduced, and the behavior
of positive solutions to the reaction-diffusion system is discussed.
Sufficient conditions for the bacteria to vanish or spread are given. We show that, if $R_0\leq 1$, the bacteria always vanish,
and if $R^F_0(t_0)\geq 1$ for some $t_0\geq 0$, the bacteria
must spread, while if $R^F_0(0)<1<R_0$, the spreading or vanishing of the bacteria depends on the initial number of  bacteria,
the length of the initial habitat, the diffusion rate, and other factors. Moreover, some sharp criteria are given. }



\end{quote}

\section{Introduction}

Recently, many mathematical models have been proposed to investigate the
spatial spread of infectious diseases epidemics  (see \cite{AM1, AM3, BHM, ZH}).
 To understand the dynamics of fecally-orally transmitted diseases in the European Mediterranean regions,
 Capasso and Maddalena  \cite{CM} have proposed an epidemic reaction-diffusion model
   described by the following coupled parabolic
system:
\begin{eqnarray}
\left\{
\begin{array}{ll}
\frac{\partial u(x,t)}{\partial t}=d\Delta u(x,t)-a_{11}u(x,t)+a_{12}v(x,t),&(x,t)\in \Omega \times (0,+\infty), \\
\frac{\partial v(x,t)}{\partial t}=-a_{22}v(x,t)+G(u(x,t)),&(x,t)\in \Omega \times (0,+\infty),\\
\frac{\partial u}{\partial \eta}+\alpha u=0, &(x,t)\in \partial\Omega \times (0,+\infty),\\
u(x,0)=u_0(x),\ v(x,0)=v_0(x),&x\in \overline \Omega,
\end{array} \right.
\label{Aa1}
\end{eqnarray}
where $u(x,t)$ and $v(x,t)$ represent the spatial densities of bacteria and the infective human population,
respectively, at a point $x$ in the habitat $\Omega\in \R^n$, and at time $t\geq 0$, and  $\partial /\partial \eta$ denotes the outward
normal derivative. The positive constant $d$ denotes
the diffusion constant of the bacteria, $1/a_{11}>0$ is the mean lifetime of the bacteria in the environment,
 the term $- a_{11}u$ denotes the natural growth rate of the bacterial population,
$1/a_{22}>0$ is the mean infectious period of an infective human,  the term $-a_{22}v$
describes the natural damping of the infective
population due to the finite mean duration of the infectiousness of humans, $a_{12}>0$ is the multiplicative factor
of the infectious bacteria due to the human population, and the term $a_{12}v$ is the
contribution of the infective humans to the growth rate of the bacteria.
The last term $G(u)$ is the infection rate of the humans under the assumption
that the total susceptible human population is constant during the evolution of the epidemic.
This kind of mechanism is used to interpret other epidemics with oro-faecal transmission
such as typhoid fever, infectious hepatitis, polyomelitis,
and the like ; see \cite{CM, Che} and the references therein for more details.

Assume that
\begin{itemize}
\item[(A1)] $G\in C^1([0, \infty)), G(0)=0, G'(z)>0, \forall z\geq 0$;

\item[(A2)] $\frac{G(z)}{z}$ is decreasing and $\lim_{z\to +\infty}\frac{G(z)}{z}<\frac {a_{11}a_{22}}{a_{12}}$.
\end{itemize}
An example is $G(z)=\frac{a_{21}z}{1+z}$ with $a_{21}>0$.

For the corresponding O.D.E. system of  (\ref{Aa1}),
\begin{eqnarray}
\left\{
\begin{array}{ll}
\frac{d u(t)}{d t}=-a_{11}u(t)+a_{12}v(t),&t>0, \\
\frac{d v(t)}{d t}=-a_{22}v(t)+G(u(t)),&t>0,
\end{array} \right.
\label{aode1}
\end{eqnarray}
linearization and spectrum analysis show that a threshold
parameter $R_0(:=\frac {G'(0)a_{12}}{a_{11}a_{22}})$ exists such that if $0 < R_0 < 1$, then the epidemic always tends to
extinction, while for $R_0 > 1$, a nontrivial endemic level appears which is globally asymptotically stable in the
positive quadrant.

For problem (\ref{Aa1}), in which the bacteria diffuse but the infective human population does not, the authors in \cite{CM}
introduced a threshold parameter $R_0^D(:=\frac {G'(0)a_{12}}{(a_{11}+d\lambda_1)a_{22}})$ such that for
$0 < R_0^D < 1$, the epidemic eventually tends to extinction, while for $R_0^D> 1$ a
globally asymptotically stable spatially inhomogeneous stationary endemic state appears, where $\lambda_1$ is
the first eigenvalue of the boundary value problem
$$-\Delta \phi=\lambda\phi \ \textrm{in}\ \Omega \ \textrm{with}\
 \frac{\partial \phi}{\partial \eta}+\alpha \phi=0\ \textrm{on}\  \partial \Omega.$$

To understand the whole dynamical structure of solutions to  (1.1) and its corresponding reaction systems, traveling waves
and entire solutions were widely studied. The existence, uniqueness and stability of traveling waves were established in
\cite{TZ, WW, WW2, WS, XZ, ZW}. Recently, Wu \cite{WS} considered entire solutions of a bistable reaction-diffusion system
 (\ref{Aa1}) in the bistable case, and  proved the existence of entire solutions that
 behave like two monotone increasing traveling wave solutions propagating
from both sides of the $x$-axis. The time-delayed and diffusive model has been considered in \cite{WW} and entire solutions have been given.
It was shown that there exist a great diversity of different types of entire solutions of reaction-diffusion equations,
which are different from traveling wave solutions.

It must be pointed out that the solution of  (\ref{Aa1}) in a fixed (bounded or unbounded) domain
is always positive for any $t>0$  no matter what the nonnegative
nontrivial initial data are. This  means that bacteria spread and the whole environment is infected immediately even
though the infection is limited to a small part of population at the beginning. This  does not match the
reality that bacteria always spread gradually. The traveling wave solutions and entire solutions play
a key role in developing a full understanding of the transient dynamics and the structure of the global attractor,
but none of those solutions can  explain the gradual expanding process.

To describe such a gradual spreading process and changing of the domain considered,
the free boundary has been introduced in many applied areas, especially
 the well-known Stefan condition used to describe the spreading process at the boundary.
 The Stefan condition was used originally
to describe the melting of ice in contact with water \cite{R} ; it was then used in
 modeling  oxygen in the muscle \cite{C},
  wound healing \cite{CF}, and more recently  the spreading of species in ecological models
 \cite{DG, DGP, DL, DL2, GW, KY, LIN, PZ,Wa}.

For emerging and re-emerging infectious bacteria,  the expanding of bacteria usually starts at
a source location and spreads over areas where
contact transmission occurs. It is crucial and interesting to study how bacteria  spread spatially
 to a larger area to cause an environmental problem.
We will focus on the changing of the infected habitat and
consider an epidemic model with the free boundary, which describes the spreading front of bacteria.  For simplicity,
assume that the human population in the whole habitat $(-\infty, \infty)$ is constant, and that the environment in $g(t)<x<h(t)$ is
infected by bacteria, the density of which is denoted by $u(x,t)$ with the infective human population denoted by $v(x,t)$,
and no bacteria or infective humans in the remaining portion of the environment.
The right spreading front of the infected environment is represented by the free boundary $x = h(t)$.
Assuming that $h(t)$ grows at a rate  proportional to the bacteria population gradient at the front \cite{LIN},
 the conditions on the right front (free boundary) are
$$u(h(t),t)=0,\quad -\mu \frac{\partial u}{\partial x}(h(t),t)=h'(t).$$
Similarly, the conditions on the left front
(free boundary) are
$$u(g(t),t)=0,\quad -\mu \frac{\partial u}{\partial x}(g(t),t)=g'(t).$$
 In such a case, we have the problem for $u(x,t)$ and $v(x,t)$ with free boundaries $x=g(t)$ and $x=h(t)$ such that
\begin{eqnarray}
\left\{
\begin{array}{ll}
\frac{\partial u(x,t)}{\partial t}=d\frac{\partial^2u(x,t)}{\partial x^2}-a_{11}u(x,t)+a_{12}v(x,t),&g(t)<x<h(t),\, t>0, \\
\frac{\partial v(x,t)}{\partial t}=-a_{22}v(x,t)+G(u(x,t)),&g(t)<x<h(t),\, t>0,\\
u(x,t)=0,\, & x=g(t)\, \textrm{or}\, x=h(t),\, t>0,\\
g(0)=-h_0,\; g'(t)=-\mu \frac{\partial u}{\partial x}(g(t), t), & t>0, \\
 h(0)=h_0, \; h'(t)=-\mu \frac{\partial u}{\partial x}(h(t), t), & t>0,\\
u(x,0)=u_0(x),\ v(x,0)=v_0(x),&-h_0\leq x\leq h_0,
\end{array} \right.
\label{a3}
\end{eqnarray}
where $x=g(t)$ and $x=h(t)$ are the moving left and right
boundaries to be determined, $h_0$ and $\mu $ are positive constants, and the initial functions
$u_{0}$ and $v_0$ are nonnegative and satisfy
\begin{eqnarray}
\left\{
\begin{array}{lll}
u_{0}\in C^2([-h_0, h_0]),\, u_{0}(\pm h_0)=0\, \textrm{and} \ 0< u_{0}(x),\, & x\in (-h_0, h_0), \\
v_{0}\in C^2([-h_0, h_0]), v_{0}(\pm h_0)=0\, \textrm{and} \ 0< v_{0}(x),\ & x\in (-h_0, h_0).
\end{array} \right.
\label{Ae1}
\end{eqnarray}

The remainder of this paper is organized as follows.
In the next section, the global existence and uniqueness of the solution to
(\ref{a3}) are proved using a contraction mapping theorem,
and a comparison principle is presented. Section 3 is devoted to sufficient conditions for the bacteria to vanish. Section 4
deals with the case and conditions for the bacteria to expand and the whole environment become infected.
Finally, we give a brief discussion in Section 5.

\section{Existence and uniqueness}

In this section, we first present the following local existence and
uniqueness result using the contraction mapping theorem and then show global existence using
suitable estimates.

\begin{thm} For any given $(u_0, v_0)$ satisfying \eqref{Ae1},  and any $\alpha \in (0, 1)$, there is a $T>0$ such that
problem  \eqref{a3} admits a unique solution
$$(u, v; g, h)\in [C^{ 1+\alpha,(1+\alpha)/2}(D_{T})]^2\times [C^{1+\alpha/2}([0,T])]^2;$$
moreover,
\begin{eqnarray} \|u\|_{C^{1+\alpha,
(1+\alpha)/2
}({D}_{T})}+\|v\|_{C^{1+\alpha,
(1+\alpha)/2}({D}_{T})}+||g\|_{C^{1+\alpha/2}([0,T])}+||h\|_{C^{1+\alpha/2}([0,T])}\leq
C,\label{b12}
\end{eqnarray}
where $D_{T}=\{(x, t)\in \R^2: x\in [g(t), h(t)], t\in [0,T]\}$, $C$
and $T$  depend only on $h_0, \alpha, \|u_0\|_{C^{2}([-h_0, h_0])}$ and $\|v_0\|_{C^{2}([-h_0, h_0])}$.
\end{thm}
\bpf
As in \cite{ZL}, we first straighten
the double free boundary fronts by making the following change of
variable:
$$y=\frac{2h_0x}{h(t)-g(t)}-\frac{h_0(h(t)+g(t))}{h(t)-g(t)},\
w(y,t)=u(x,t),\ z(y,t)=v(x,t).$$ \\
Then  (\ref{a3}) can be transformed into
\begin{eqnarray}
\left\{
\begin{array}{lll}
w_{t}=Aw_{y}+Bw_{yy}-a_{11}w(y,t)+a_{12}z(y,t),\; &t>0, \ -h_0<y<h_0, \\
z_{t}=Az_y-a_{22}z(y,t)+G(w(y,t)),\; &t>0, \ -h_0<y<h_0, \\
w=0,\quad h'(t)=-\frac{2h_0\mu}{h(t)-g(t)}\frac{\partial w}{\partial
y},\quad &t>0, \ y=h_0,\\
w=0,\quad g'(t)=-\frac{2h_0\mu}{h(t)-g(t)}\frac{\partial w}{\partial
y},\quad &t>0, \ y=-h_0,\\
h(0)=h_0, \quad g(0)=-h_0, &\\
w(y,0)=w_0(y):=u_{0}(y),\; z(y,0)=z_0(y):=v_{0}(y),\; &-h_0\leq y\leq
h_0,
\end{array} \right.
\label{Hb}
\end{eqnarray}
where $A=A(h, g,
y)=y\frac{h'(t)-g'(t)}{h(t)-g(t)}+h_0\frac{h'(t)+g'(t)}{h(t)-g(t)}$, and
$B=B(h, g)=\frac{4h_0^2d}{(h(t)-g(t))^2}$. This transformation
changes the free boundaries $x=h(t)$ and $x=g(t)$ to the fixed lines
$y=h_0$ and $y=-h_0$ respectively; therefore, the equations become more
complex, because now the coefficients in the first and second equations of
(\ref{Hb}) contain unknown functions $h(t)$ and $g(t)$.

 The rest of the  proof uses by the contraction mapping argument as in \cite{DL, ZL}
 with suitable modifications; we omit it here.
\epf

To show the global existence of the solution, we need the following estimate.

\begin{lem} Let $(u, v; g, h)$ be a solution to  \eqref{a3} defined for $t\in (0,T_0]$ for some $T_0\in (0, +\infty)$.
Then there exist constants $C_1$ and $C_2$ independent of $T_0$ such
that
\[
0<u(x, t)\leq C_1\; \mbox{ for }\, g(t)<x<h(t),\; t\in (0, T_0], \]
\[
0<v(x, t)\leq C_2\; \mbox{ for }\, g(t)<x<h(t),\; t\in (0, T_0]. \]
\end{lem}
\bpf The positivity of $u$ and $v$ are obvious, since  the initial values are nontrivial and nonnegative and the system is
quasi-increasing. Now let us consider its upper bounds. Note that $\lim_{z\to +\infty}\frac{G(z)}{z}<\frac {a_{11}a_{22}}{a_{12}}$
by the assumption $(A2)$; therefore,
there exist $C_1$ and $C_2$ such that
$$C_1\geq u_0(x),\ C_2\geq v_0(x)\, \textrm{in}\, [-h_0, h_0],$$
$$-a_{11}C_1+a_{12}C_2<0,\ -a_{22}C_2+G(C_1)<0.$$
Define
$$M_1=M_1(u, T_0)=\max_{[g(t),h(t)]\times [0, T_0]} u(x,t),$$
$$M_2=M_2(G, C_1, T_0)=\max_{0\leq z\leq \max\{C_1, M_1\}} G'(z)$$
and let
$(U(x,t), V(x,t))=(C_1-u, C_2-v)e^{-(a_{12}+M_2)t}$,  then we have
\begin{eqnarray*}
G(u)&=&G(C_1-Ue^{(a_{12}+M_2)t})\\
&=&G(C_1)-G'(\xi(x,t))Ue^{(a_{12}+M_2)t},
\end{eqnarray*}
where $\xi(x,t)$ is between $C_1$ and $u(x,t)$, and therefore $(U, V)$ satisfies
\begin{eqnarray}
\left\{
\begin{array}{ll}
\frac{\partial U (x,t)}{\partial t}> d\frac{\partial^2U(x,t)}{\partial x^2}-(a_{11}+a_{12}+M_2)U(x,t)&\\
\qquad +a_{12}V(x,t),&g(t)<x<h(t),\, 0<t\leq T_0, \\
\frac{\partial V(x,t)}{\partial t}>-(a_{22}+a_{12}+M_2)V(x,t)&\\
\qquad +G'(\xi(x,t))U(x,t),&g(t)<x<h(t),\, 0<t\leq T_0,\\
U(x,t)=C_1e^{-(a_{12}+M_2)t}, & x=g(t)\, \textrm{or}\, x=h(t),\, t>0,\\
V(x, t)=C_2e^{-(a_{12}+M_2)t}, & x=g(t)\, \textrm{or}\, x=h(t),\, t>0,\\
U(x,0)\geq 0,\ V(x,0)\geq 0,&-h_0\leq x\leq h_0.
\end{array} \right.
\label{eqf}
\end{eqnarray}

We now show that $\min \{U(x,t),V(x,t)\}\geq 0$ in $[g(t), h(t)]\times [0, T_0]$. Otherwise, there exists $(x_0, t_0)\in
(g(t), h(t))\times (0, T_0]$ such that
$$\min \{U(x_0,t_0),V(x_0,t_0)\}=\min_{[g(t), h(t)]\times [0, T_0]} \min \{U(x,t),V(x,t)\}<0.$$

If $U(x_0, t_0)=\min \{U(x_0,t_0),V(x_0,t_0)\}<0$, then $U(x,t)$ attains its minimum in $[g(t), h(t)]\times [0, T_0]$ at $(x_0,t_0)$;
therefore,
$$\frac{\partial U (x_0,t_0)}{\partial t}- d\frac{\partial^2U(x_0,t_0)}{\partial x^2}\leq 0.$$
However,  $-(a_{11}+a_{12}+M_2)U(x_0,t_0)+a_{12}V(x_0,t_0)\geq -(a_{11}+M_2)U(x_0,t_0)>0$, which leads a contradiction to the first inequality in (\ref{eqf}).

Similarly, if $V(x_0, t_0)=\min \{U(x_0,t_0),V(x_0,t_0)\}<0$, then $V(x,t)$ attains its minimum in $[g(t), h(t)]\times [0, T_0]$ at $(x_0,t_0)$;
therefore, $\frac{\partial V (x_0,t_0)}{\partial t}\leq 0$.
However, $-(a_{22}+a_{12}+M_2)V(x_0,t_0)+G'(\xi)U(x_0,t_0)\geq -(a_{22}+a_{12})V(x_0,t_0)>0$, which leads to a contradiction to the second inequality in (\ref{eqf}).
Thus, we have $\min \{U(x,t),V(x,t)\}\geq 0$ in $[g(t), h(t)]\times [0, T_0]$, or $u(x,t)\leq C_1$ and $v(x,t)\leq C_2$ in $[g(t), h(t)]\times [0, T_0]$.

\epf

The next lemma shows that the left free boundary for (\ref{a3}) is strictly monotone decreasing and the right boundary is increasing.

\begin{lem} Let $(u, v; g, h)$ be a solution to \eqref{a3} defined for $t\in (0,T_0]$ for some $T_0\in (0, +\infty)$.
Then there exists a constant $C_3$  independent of $T_0$ such
that
\[
 0<-g'(t),\ h'(t)\leq C_3 \; \mbox{ for } \; t\in (0,T_0]. \]
\end{lem}
\bpf
Applying the strong maximum principle to the equation of $u$ gives
 $$u_x(h(t), t)<0 \ \;\; \textrm{for} \ 0<t\leq T_0.$$
 Hence $h'(t)>0$ for $t\in (0, T_0]$ by the free boundary condition in (\ref{a3}). Similarly, $g'(t)<0$ for $t\in (0, T_0]$.

It remains to be shown that $-g'(t),\ h'(t)\leq C_3$ for $t\in (0,T_0]$ and
some $C_3$. The proof is similar to that of Lemma 2.2 in \cite{DL} with $C_3=2MC_1\mu$ and
$$M=\max\left\{ \frac 1{h_0},\ \sqrt{\frac{a_{12}C_2}{2dC_1}},\
\frac{4\|u_0\|_{C^1([-h_0,h_0])}}{3C_1}\right\},$$ we omit it here.
 \epf

Since  $u,v$ and $g'(t), h'(t)$ are bounded in $(g(t),h(t))\times (0, T_0]$ by constants independent of $T_0$, the global solution is guaranteed.

\begin{thm} The solution of  \eqref{a3} exists and is
unique for all $t\in (0,\infty)$.
\end{thm}

In what follows, we exhibit the comparison principle, which can be proved similarly to a Lemma 3.5 in \cite{DL}.
\begin{lem} (The Comparison Principle)
  Assume that $\overline g, \overline
h\in C^1([0, +\infty))$, $\overline u(x, t), \overline v(x,t)\in  C([\overline g(t), \overline h(t)]\times [0, +\infty))\cap
C^{2,1}((\overline g(t), \overline h(t))\times (0, +\infty))$, and
\begin{eqnarray*}
\left\{
\begin{array}{lll}
\frac{\partial \overline u}{\partial t}\geq d\frac{\partial^2 \overline u}{\partial x^2}-a_{11}\overline u+a_{12}
\overline v,&\overline g(t)<x<\overline h(t),\, t>0, \\
\frac{\partial \overline v}{\partial t}\geq -a_{22}\overline v+G(\overline u),&\overline g(t)<x<\overline h(t),\, t>0,\\
\overline u(x,t)=\overline v(x, t)=0,\, & x=\overline g(t)\, \textrm{or}\, x=\overline h(t),\, t>0,\\
\overline g(0)\leq -h_0,\; \overline g'(t)\leq -\mu \frac{\partial \overline u}{\partial x}(\overline g(t), t), & t>0, \\
 \overline h(0)\geq h_0, \; \overline h'(t)\geq -\mu \frac{\partial \overline u}{\partial x}(\overline h(t), t), & t>0,\\
\overline u(x,0)\geq u_0(x),\ \overline v(x,0)\geq v_0(x),&-h_0\leq x\leq h_0.
\end{array} \right.
\end{eqnarray*}
Then the solution $(u, v; g, h)$ to the free boundary problem $(\ref{a3})$ satisfies
$$h(t)\leq\overline h(t),\ g(t)\geq \overline g(t),\quad t\in [0, +\infty),$$
$$u(x, t)\leq \overline u(x, t),\ v(x, t)\leq \overline v(x, t),\quad (x, t)\in [g(t), h(t)]\times [0, +\infty).$$\label{Com}
\end{lem}
\begin{rmk} The pair $(\overline u, \overline h)$ in Lemma \ref{Com} is usually called an upper solution
of  \eqref{a3}. We can define a lower solution by
reversing all of the inequalities in the obvious places. Moreover, one
can easily prove an analogue of Lemma \ref{Com} for lower solutions.
\end{rmk}

We next fix $v_0, \mu, a_{ij}$, let $u_0=\sigma \phi(x)$ and examine the dependence of the solution on $\sigma$,
 writing $(u^{\sigma}, v^{\sigma}; g^{\sigma}, h^{\sigma})$ to emphasize this dependence. As a corollary of Lemma \ref{Com},
we have the following monotonicity:

\begin{cor} Let $(u_0, v_0)=\sigma (\phi(x),\psi(x))$. For fixed $\phi(x), \psi(x), \mu$ and $a_{ij}$,
if $\sigma_1\leq \sigma_2$, then $u^{\sigma_1}(x, t)\leq u^{\sigma_2}(x, t)$ and
$v^{\sigma_1}(x, t)\leq v^{\sigma_2}(x, t)$ in $[g^{\sigma_1}(t), h^{\sigma_1}(t)]\times (0, \infty)$,
$g^{\sigma_1}(t)\geq g^{\sigma_2}(t)$ and $h^{\sigma_1}(t)\leq h^{\sigma_2}(t)$ in $(0, \infty)$.
\end{cor}

\section{Bacteria vanishing}
It follows from Lemma 2.3 that $x=h(t)$ is monotonic increasing, $x=g(t)$ is monotonic decreasing and therefore
there exist  $h_\infty, -g_\infty\in (0, +\infty]$ such that $\lim_{t\to +\infty} \ h(t)=h_\infty$
and $\lim_{t\to +\infty} \ g(t)=g_\infty$. The next lemma shows that if $h_\infty<\infty$, then $-g_\infty<\infty$,
and vice versa. That is, the double free boundary
fronts $x=g(t)$ and $x=h(t)$ are both finite or infinite simultaneously.
\begin{lem} Let $(u, v; g, h)$ be a solution to  $(\ref{a3})$
defined for $t\in[0, +\infty)$ and $x\in[g(t), h(t)]$. Then we have
$$-2h_0<g(t)+h(t)<2h_0 \mbox{ for } t\in[0, +\infty).$$
\end{lem}
\bpf By continuity we know $g(t)+h(t)>-2h_0$ holds for small $t>0$.
Define
$$T:=\sup\{s: g(t)+h(t)>-2h_0 \mbox{ for all }  t\in[0,s)\}.$$ As in \cite{DB}, we claim
that $T=\infty$. Otherwise, $0<T<\infty$  and
$$g(t)+h(t)>-2h_0 \mbox{ for } t\in[0,T),\ g(T)+h(T)=-2h_0.$$
Hence,
\begin{eqnarray}
g'(T)+h'(T)\leq0. \label{Hq}
\end{eqnarray}

To get a contradiction, we consider the functions $$w(x, t):=u(x,t)-u(-x-2h_0, t),\ z(x, t):=v(x, t)-v(-x-2h_0, t)$$
over the region $$\Lambda:=\{(x, t):\ x\in[g(t), -h_0],\, t\in[0, T]\}.$$
It is easy  to check that the pair $(w, z)$ is
well-defined for $(x, t)\in\Lambda$ since $-h_0\leq -x-2h_0\leq -g(t)-2h_0\leq h(t)$, and the pair satisfies
$$w_t-d w_{xx}=-a_{11} w+a_{12}z \mbox{ for } g(t)<x<-h_0,\ 0<t\leq T,$$
$$z_t=c_{21}(x,t)w-a_{22} z\mbox{ for } g(t)<x<-h_0,\ 0<t\leq T$$
with $0\leq c_{21}:= \frac {G(u(x,t))-G(u(-x-2h_0,t))}{u(x,t)-u(-x-2h_0,t)}\in L^\infty(\Lambda)$, and
$$w(-h_0, t)=z(-h_0,t)=0,\, w(g(t), t)< 0,\, z(g(t), t)< 0 \mbox{ for } 0<t<T.$$ Moreover,
$$w(g(T), T)=u(g(T), T)-u(-g(T)-2h_0, T)=u(g(T), T)-u(h(T), T)=0.$$
Applying the proof for the strong maximum principle and the Hopf lemma, we deduce
$$w(x,t)<0,\, z(x,t)<0  \mbox{ in } (g(t),-h_0)\times (0, T] \mbox{ and } w_x(g(T), T)<0.$$
However, $$w_x(g(T), T)=\frac {\partial u}{\partial x}(g(T), T)+\frac {\partial u}{\partial x}(
h(T), T)=-[g'(T)+h'(T)]/(\mu),$$ which implies $$g'(T)+h'(T)>0,$$ a
contradiction to (\ref{Hq}). Hence we have proven $$g(t)+h(t)>-2h_0 \mbox{
for all } t>0.$$

Analogously, we can prove $g(t)+h(t)<2h_0 \mbox{ for all } t>0$ by considering
$$W(x, t):=u(x,t)-u(2h_0-x, t),\ Z(x, t):=v(x, t)-v(2h_0-x, t)$$ over the
region $\Lambda':=[h_0, h(t)]\times [0, T']$ with
$T':=\sup\{s: g(t)+h(t)<2h_0 \mbox{ for all }  t\in[0,s)\}$. The
completes the proof.  \epf

 Next, we discuss the properties of the free boundary, because
 the transmission of the bacteria depends on whether  $h_\infty-g_\infty=\infty$ and $\limsup_{t\to
+\infty}\ (||u(\cdot, t)||_{C(g(t), h(t)])}+||v(\cdot, t)||_{C([g(t),h(t)])})=0$. We then have the following definitions:

\begin{defi}
The bacteria are {\bf vanishing} if
$$h_\infty-g_\infty <\infty\ \textrm{ and}\
 \lim_{t\to +\infty} \ (||u(\cdot, t)||_{C([g(t),h(t)])}+||v(\cdot, t)||_{C([g(t), h(t)])})=0,$$
  and  {\bf spreading} if $$h_\infty-g_\infty =\infty\ \textrm{and}\
\limsup_{t\to +\infty}\ (||u(\cdot, t)||_{C([g(t),h(t)])}+||v(\cdot, t)||_{C([g(t),h(t)])})>0.$$
\end{defi}

 The next result shows that if $h_\infty-g_\infty<\infty$, then vanishing occurs.
\begin{lem}  If $h_\infty-g_\infty<\infty$, then $\lim_{t\to
+\infty} \ (||u(\cdot, t)||_{C([g(t),h(t)])}+||v(\cdot, t)||_{C([g(t),h(t)])})=0$.
\end{lem}
\bpf We first prove that $\lim_{t\to
+\infty} \ ||u(\cdot, t)||_{C([g(t),h(t)])}=0$. Assume that
$$\limsup_{t\to +\infty} \ ||u(\cdot, t)||_{C([g(t), h(t)])}=\delta>0$$
 by contradiction. Then there exists a sequence $(x_k, t_k )$
in $(g(t), h(t))\times (0, \infty)$
such that $u(x_k,t_k)\geq \delta /2$ for all $k \in \mathbb{N}$, and $t_k\to \infty$ as $k\to \infty$.
Since  $-\infty<g_\infty<g(t)<x_k<h(t)<h_\infty<\infty$, we then have that a subsequence of $\{x_n\}$ converges
to $x_0\in (g_\infty, h_\infty)$. Without loss of generality, we assume $x_k\to x_0$ as $k\to \infty$.

Define $W_k(x,t)=u(x,t_k+t)$ and $Z_k(x,t)=v(x,t_k+t)$  for
$x\in (g(t_k+t), h(t_k+t)), t\in (-t_k, \infty)$.
It follows from  parabolic regularity that  $\{(W_k, Z_k)\}$ has a subsequence $\{(W_{k_i}, Z_{k_i})\}$ such that
$(W_{k_i}, Z_{k_i})\to (\tilde W, \tilde Z)$ as $i\to \infty$ and $(\tilde W, \tilde Z)$ satisfies
\begin{eqnarray*} \left\{
\begin{array}{lll}
\tilde W_t-d \tilde W_{xx}=-a_{11}\tilde W+a_{12}\tilde Z,\; & g_\infty<x<h_\infty,\ t\in (-\infty, \infty),  \\
\tilde Z_t=-a_{22}\tilde Z+G(\tilde W),\; &\ g_\infty<x<h_\infty, \ t\in (-\infty, \infty).
\end{array} \right.
\end{eqnarray*}
Note that $\tilde W(x_0, 0)\geq \delta/2$; therefore, $\tilde W>0$ in $ (g_\infty, h_\infty)\times(-\infty, \infty)$.

 Using a similar method to prove the Hopf lemma at the point $(h_\infty, 0)$ yields
 $\tilde W_x(h_\infty, 0 )\leq -\sigma_0$ for some $\sigma_0>0$.

On the other hand, since $-g(t)$ and $h(t)$ are increasing and bounded, it follows from standard $L^p$ theory and then the Sobolev imbedding
theorem (\cite{LSU}) that for any $0<\alpha <1$, there exists a constant $\tilde C$
depending on $\alpha, h_0, \|u_{0}\|_{C^{2}[-h_0, h_0]}$, $\|v_{0}\|_{C^{2}[-h_0, h_0]}$, and $g_\infty, h_\infty$ such that
\begin{eqnarray}\|u\|_{C^{1+\alpha,
(1+\alpha)/2}([g(t), h(t)]\times[0, \infty))}+\|h\|_{C^{1+\alpha/2}([0,\infty))}\leq
\tilde C.\label{Bg1}
\end{eqnarray}

Now, since $\|h\|_{C^{1+\alpha/2}([0,\infty))}\leq
\tilde C$ and $h(t)$ is bounded,  we then have $h'(t)\to 0$ as $t\to \infty$, that is,
$\frac {\partial u}{\partial x}(h(t_k),t_k)\to 0$ as $t_k\to \infty$ by the free boundary condition. Moreover,
the fact that  $\|u\|_{C^{
1+\alpha,(1+\alpha)/2}([g(t), h(t)]\times[0, \infty))}\leq \tilde C$ gives
$\frac {\partial u}{\partial x}(h(t_k),t_k+0)=(W_k)_x(h(t_k),0)\to \tilde W_x(h_\infty,0)$ as $k\to \infty$, and then $\tilde W_x(h_\infty,0)=0$,
which leads to a contradiction to the fact that  $\tilde W_x(h_\infty,0)\leq -\sigma_0<0$.
 Thus $\lim_{t\to +\infty} \ ||u(\cdot,t)||_{C([g(t),h(t)])}=0$.

 Note that $v(x,t)$ satisfies
 $$\frac{\partial v(x,t)}{\partial t}=-a_{22}v(x,t)+G(u(x,t)),\ g(t)<x<h(t),\, t>0,$$
 and $G(u(x,t))\to 0$ uniformly for $x\in [g(t), h(t)]$ as $t\to \infty$; therefore, we have $\lim_{t\to
+\infty} \ ||v(\cdot, t)||_{C([g(t),h(t)])}=0$.
\epf

In the introduction, a threshold $R_0$, usually called the basic reproduction number, is given to decide whether the bacteria described by  (\ref{aode1}) vanish.
Notice that the interval domain for free boundary problem \eqref{a3} changes with $t$; therefore,
the basic reproduction number is not a constant and
should change with $t$.

Now we introduce the basic reproduction number $R_0^F(t)$ for  \eqref{a3} by
$$R_0^F(t):=R_0^D((g(t), h(t)))=\frac {G'(0)\frac {a_{12}}{a_{22}}}{a_{11}+d(\frac \pi {h(t)-g(t)})^2},$$
where we use $R_0^D(\Omega)$ to denote the basic reproduction number for the corresponding problem in $\Omega$
with null Dirichlet boundary condition on $\partial \Omega$.
Now,  the following result is obvious; see also Lemma 2.3 in \cite{HH}.
\begin{lem} $1-R_0^F(t)$ has the same sign as $\lambda_0$, where $\lambda_0$ is the principal eigenvalue of the problem
\begin{eqnarray}
\left\{
\begin{array}{lll}
-d \psi_{xx}=-a_{11} \psi+{G'(0)\frac {a_{12}}{a_{22}}}\psi+\lambda_0 \psi,\; &
x\in (g(t), h(t)),  \\
\psi(x)=0, &x=g(t)\ \textrm{or}\ x=h(t).
\end{array} \right.
\label{B1f}
\end{eqnarray}
\end{lem}

In fact, here
$$\lambda_0=a_{11}+d(\frac \pi {h(t)-g(t)})^2-G'(0)\frac {a_{12}}{a_{22}}=[a_{11}+d(\frac \pi {h(t)-g(t)})^2](1-R_0^F(t)).$$
With the above defined reproduction number, we also have
\begin{lem} $R_0^F(t)$ is strictly monotone increasing function of $t$, that is if $t_1<t_2$, then $R_0^F(t_1)<R_0^F(t_2)$.
Moreover, if $h(t)\to \infty$ as $t\to \infty$, then $R_0^F(t)\to R_0$ as $t\to \infty$.
\end{lem}

Next we give sufficient conditions so that the bacteria are vanishing.
\begin{thm} If $R_0\leq 1$, then $h_\infty-g_\infty<\infty$ and
$\lim_{t\to +\infty} \ (||u(\cdot, t)||_{C([g(t),h(t)])}+||v(\cdot, t)||_{C([g(t),h(t)])})=0$.
\end{thm}
\bpf
We first show that $h_\infty -g_\infty<+\infty$. In fact, direct calculations yield
\begin{eqnarray*}& &\frac{\textrm{d}}{\textrm{d} t}\int_{g(t)}^{h(t)}[u(x, t)+\frac {a_{12}}{a_{22}}v(x, t)]\textrm{d}x\\[1mm]
&=&\int_{g(t)}^{h(t)}[u_t+\frac {a_{12}}{a_{22}}v_t](x, t)\textrm{d}x+h'(t)[u+\frac {a_{12}}{a_{22}}v](h(t), t)
-g'(t)[u+\frac {a_{12}}{a_{22}}v](g(t), t)\\[1mm]
&=&\int_{g(t)}^{h(t)}d u_{xx}\textrm{d}x+\int_{g(t)}^{h(t)}-a_{11}u(x,t)+\frac {a_{12}}{a_{22}}G(u(x, t))\textrm{d}x\\[1mm]
&=&-\frac{d}{\mu}(h'(t)-g'(t))+\int_{g(t)}^{h(t)}-a_{11}u(x,t)+\frac {a_{12}}{a_{22}}G(u(x, t))\textrm{d}x.
\end{eqnarray*}
Integrating from $0$ to $t\,(>0)$ gives
\begin{eqnarray}
& &\int_{g(t)}^{h(t)}[u+\frac {a_{12}}{a_{22}}v](x, t)\textrm{d}x = \int ^{h(0)}_{g(0)}[u+\frac {a_{12}}{a_{22}}v](x, 0)\textrm{d}x \\
& &\quad +\frac {d}{\mu}(h(0)-g(0))-\frac {d}{\mu}(h(t)-g(t))\nonumber \\[1mm]
& &\quad +\int_{0}^t\int_{g(s)}^{h(s)}-a_{11}u(x,t)+\frac {a_{12}}{a_{22}}G(u(x, t))dxds,
\quad t\geq 0.\label{k1}
\end{eqnarray}
Since $\frac {G(z)}{z}\leq G'(0)$ by the assumption $(A2)$, it  follows from $R_0\leq 1$ that $-a_{11}u(x,t)+\frac {a_{12}}{a_{22}}G(u(x, t))\leq 0$ for $x\in [g(t),h(t)]$ and $t\geq 0$, we have
$$\frac {d}{\mu}(h(t)-g(t)) \leq \int ^{h(0)}_{g(0)}[u+\frac {a_{12}}{a_{22}}v](x, 0)\textrm{d}x
+\frac {d}{\mu}(h(0)-g(0))$$
for $t\geq 0$, which in turn gives that $h_\infty-g_\infty<\infty$. Therefore, the bacteria are vanishing
 as a consequence of  Lemma 3.2.
\epf

\begin{thm} If $R_0^F(0)<1$ and $||u_0(x)||_{C([-h_0, h_0])}$, $||v_0(x)||_{C([-h_0, h_0])}$ are sufficiently small. Then $h_\infty-g_\infty<\infty$ and
$\lim_{t\to +\infty} \ (||u(\cdot, t)||_{C([g(t),h(t)])}+||v(\cdot, t)||_{C([g(t),h(t)])})=0$.
\end{thm}
\bpf We construct a suitable upper solution for $(u, v)$.
Since $R_0^F(0)<1$, it follows from Lemma 3.3 that there is a $\lambda_0>0$ and $0<\psi(x)\leq 1$ in $(-h_0, h_0)$ such that
\begin{eqnarray}
\left\{
\begin{array}{lll}
-d \psi_{xx}= -a_{11}\psi+{G'(0)\frac {a_{12}}{a_{22}}}\psi+\lambda_0 \psi,\; &
-h_0<x<h_0,  \\
\psi(x)=0, &x=\pm h_0.
\end{array} \right.
\label{B1f1}
\end{eqnarray}
Therefore, there exists a small $\delta >0$ such that
$$-\delta +(\frac 1{(1+\delta)^2}-1)|-a_{11}+{G'(0)\frac {a_{12}}{a_{22}}}|+[\frac 1{(1+\delta)^2}-\frac 14]\lambda_0\geq 0.$$

Similarly as in \cite{DL}, we set
$$\sigma (t)=h_0(1+\delta-\frac \delta 2 e^{-\delta t}), \  t\geq 0,$$
and
$$\overline u(x, t)=\varepsilon e^{-\delta t}\psi(xh_0/\sigma (t)), \ -\sigma(t)\leq
x\leq \sigma(t),\ t\geq 0.$$
$$\overline v(x,t)=(\frac {G'(0)}{a_{22}}+\frac{\lambda_0}{4a_{12}})\overline u(x,t), \ -\sigma(t)\leq
x\leq \sigma(t),\ t\geq 0.$$

Direct computations yield
\begin{eqnarray*}
& &\overline u_t-d \dfrac{\partial^2\overline u}{\partial x^2}+
a_{11}\overline u-a_{12}\overline v\\
& &=-\delta \overline u-\varepsilon e^{-\delta t}\psi'\frac{xh_0\sigma'(t)}{\sigma^2(t)}+(\frac{h_0}{\sigma(t)})^2
[-a_{11}+{G'(0)\frac {a_{12}}{a_{22}}}+\lambda_0 ]\overline u\\
& & \quad +[a_{11}-{G'(0)\frac {a_{12}}{a_{22}}}-\frac{\lambda_0}{4}]\overline u\\
& &\geq \overline u \{-\delta +(\frac 1{(1+\delta)^2}-1)|-a_{11}+{G'(0)\frac {a_{12}}{a_{22}}}|+[\frac 1{(1+\delta)^2}-\frac 14]\lambda_0\}\geq 0,
\end{eqnarray*}
\begin{eqnarray*}
& &\overline v_t+a_{22}\overline v-G(\overline u(x,t))\\
& &=-\delta \overline v-\varepsilon e^{-\delta t}\psi'\frac{xh_0\sigma'(t)}{\sigma^2(t)}(\frac {G'(0)}{a_{22}}+\frac{\lambda_0}{4a_{12}})
+a_{22}(\frac {G'(0)}{a_{22}}+\frac{\lambda_0}{4a_{12}})\overline u(x,t)-G(\overline u(x,t))\\
& &\geq (a_{22}-\delta )(\frac {G'(0)}{a_{22}}+\frac{\lambda_0}{4a_{12}})\overline u(x,t)-G(\overline u(x,t))\\
& &=(a_{22}-\delta )(\frac {G'(0)}{a_{22}}+\frac{\lambda_0}{4a_{12}})\overline u(x,t)-G'(\xi(x,t))\overline u(x,t)\\
& &=[G'(0)-G'(\xi(x,t))+ (a_{22}-\delta )\frac{\lambda_0}{4a_{12}}-\frac {G'(0)}{a_{22}}\delta]\overline u(x,t)
\end{eqnarray*}
for all $t>0$ and $-\sigma (t)<x<\sigma (t)$, where $\xi\in (0, \overline u)$.  Since $\overline u\leq \varepsilon$,
 if
$\delta$ and $\varepsilon$ are sufficiently small, then we have
$$[G'(0)-G'(\xi(x,t))+ (a_{22}-\delta )\frac{\lambda_0}{4a_{12}}-\frac {G'(0)}{a_{22}}\delta]\geq 0.$$
On the other hand, we have
$\sigma'(t)=h_0 \frac {\delta^2} 2 e^{-\delta t}$, $-\overline u_x(
\sigma (t),t)=-\varepsilon \frac {h_0}{\sigma (t)}\psi'(h_0)e^{-\delta t}$, and $-\overline u_x(
-\sigma (t),t)=-\varepsilon \frac {h_0}{\sigma (t)}\psi'(-h_0)e^{-\delta t}$. Noticing that $\psi'(-h_0)=-\psi'(h_0)$,
we now choose $\varepsilon=-\frac {\delta^2h_0(1+\delta)} {2\mu \psi'(h_0)}$ such that
\begin{eqnarray*}
\left\{
\begin{array}{lll}
\dfrac{\partial \overline u}{\partial t}\geq d \dfrac{\partial^2 \overline u}{\partial x^2}
-a_{11}\overline u+a_{12}\overline v,\; & -\sigma(t)<x<\sigma(t),\, t>0, \\
\dfrac{\partial \overline v}{\partial t}\geq -a_{22}\overline v+G(\overline u),\; &  -\sigma(t)<x<\sigma(t), \, t>0,\\
\overline u(x,t)=\overline v(x, t)=0,&x=\pm \sigma(t)\, \, t>0,\\
-\sigma (0)<-h_0,\; -\sigma'(t)\leq -\mu \frac{\partial \overline u}{\partial x}(-\sigma(t), t), & t>0, \\
\sigma(0)> h_0, \; \sigma'(t)\geq -\mu \frac{\partial \overline u}{\partial x}(\sigma(t), t), & t>0.
\end{array} \right.
\end{eqnarray*}
If  $||u_{0}||_{L^\infty}\leq \varepsilon \psi(\frac {h_0}{1+\delta/2})$
and $||v_{0}||_{L^\infty}\leq \varepsilon \psi(\frac {h_0}{1+\delta/2})(\frac {G'(0)}{a_{22}}+\frac{\lambda_0}{4a_{12}})$, then
$u_{0}(x)\leq \varepsilon \psi(\frac {h_0}{1+\delta/2})\leq \overline u(x, 0)=\varepsilon \psi(\frac {x}{1+\delta/2})$
and $v_{0}(x)\leq \varepsilon \psi(\frac {h_0}{1+\delta/2})(\frac {G'(0)}{a_{22}}+\frac{\lambda_0}{4a_{12}})\leq \overline v(x, 0)$
 for $x\in [-h_0, h_0]$.
We can now apply Lemma 2.5 to conclude that $g(t)\geq -\sigma(t)$ and $h(t)\leq\sigma(t)$ for $t>0$. It
follows that $h_\infty-g_\infty\leq \lim_{t\to\infty}
2\sigma(t)=2h_0(1+\delta)<\infty$, and $\lim_{t\to +\infty} \ (||u(\cdot, t)||_{C([g(t),h(t)])}+||v(\cdot, t)||_{C([g(t),h(t)])})=0$ by Lemma 3.2.
 \epf

\section{Bacteria spreading}

In this section, we  give the sufficient conditions for the bacteria to be spreading. We first
prove that if $R_0^F(0)\geq 1$, the bacteria are spreading.
\begin{thm} If $R_0^F(0)\geq 1$, then $h_\infty=-g_\infty=\infty$ and $\liminf_{t\to
+\infty} \ ||u(\cdot, t)||_{C([0, h(t)])}>0$, that is, spreading occurs.
\end{thm}
\bpf We first consider the case that $R_0^F(0):=R_0^D((-h_0, h_0))>1$. In this case, we have that the eigenvalue problem
\begin{eqnarray}
\left\{
\begin{array}{lll}
-d \psi_{xx}= -a_{11}\psi+{G'(0)\frac {a_{12}}{a_{22}}}\psi+\lambda_0 \psi,\; &
-h_0<x<h_0,  \\
\psi(x)=0, &x=\pm h_0
\end{array} \right.
\label{B2f}
\end{eqnarray}
admits a positive solution $\psi(x)$ with $||\psi||_{L^\infty}=1$, where $\lambda_0$ is the principal eigenvalue. It follows from Lemma 3.3
that $\lambda_0<0$.

We  construct a suitable lower solution to
\eqref{a3}, and we define
$$\underline {u}(x,t)=\delta \psi(x),\quad \underline v=(\frac {G'(0)}{a_{22}}+\frac{\lambda_0}{4a_{12}})\delta \psi(x)$$
for $-h_0\leq x\leq h_0$, $t\geq 0$, where $\delta $ is sufficiently small.

 Direct computations yield
\begin{eqnarray*}
& &\dfrac{\partial \underline u}{\partial t}-d \dfrac{\partial^2\underline u}{\partial x^2}+a_{11}
\underline u-a_{12}\underline v=\delta \psi(x)(\frac 34\lambda_0)\leq 0\\
& &\dfrac{\partial \underline v}{\partial t}+a_{22}\underline v-G(\underline u) =\delta \psi(x)[G'(0)-G'(\xi(x,t))+\frac {a_{22}\lambda_0}{4a_{12}}]
\end{eqnarray*}
for all $t>0$ and $-h_0<x<h_0$, where $\xi\in (0, \underline u)$. Noting that $\lambda_0<0$ and $0\leq \xi(x,t)\leq \underline u(x,t)\leq \delta$, we can chose $\delta $ sufficiently small such that
\begin{eqnarray*}
\left\{
\begin{array}{lll}
\dfrac{\partial \underline u}{\partial t}\leq d \dfrac{\partial^2 \underline u}{\partial x^2}
-a_{11}\underline u+a_{12}\underline v,\; & -h_0<x<h_0,\, t>0, \\
\dfrac{\partial \underline v}{\partial t}\leq -a_{22}\underline v+G(\underline u),\; &  -h_0<x<h_0, \, t>0,\\
\underline u(x,t)=\underline v(x, t)=0,&x=\pm h_0\, \, t>0,\\
0=-h'_0\geq -\mu \frac{\partial \underline u}{\partial x}(-h_0, t), & t>0, \\
0=h'_0\leq -\mu \frac{\partial \underline u}{\partial x}(h_0, t), & t>0,\\
\underline {u}(x,0)\leq u_{0}(x),\ \underline{v}(x,0)\leq v_{0}(x),\; &-h_0\leq x\leq h_0.
\end{array} \right.
\end{eqnarray*}
Hence, applying Remark 2.1 yields that $u(x,t)\geq\underline u(x,t)$ and  $v(x,t)\geq\underline v(x,t)$
in $[-h_0, h_0]\times [0,\infty)$. It follows that $\liminf_{t\to
+\infty} \ ||u(\cdot, t)||_{C([g(t), h(t)])}\geq \delta \psi(0)>0$ and therefore $h_\infty-g_\infty=+\infty$ by Lemma 3.2.

If $R_0^F(0)=1$,  then for any positive time $t_0$, we have $g(t_0)<-h_0$ and $h(t_0)>h_0$; therefore, $R_0^F(t_0)>
R_0^F(0)=1$ by the monotonicity in Lemma 3.4.
Replacing the initial time $0$ by the positive time $t_0$, we then have $h_\infty-g_\infty=+\infty$ as above.
 \epf

\begin{rmk} It follows from the above proof that spreading occurs, if there exists $t_0\geq 0$ such that $R_0^F(t_0)\geq 1$.
\end{rmk}

Theorem 3.6 shows if $R_0^F(0)<1$, vanishing occurs for small initial size of infected bacteria, and Theorem 3.5 implies that
if $R_0\leq 1$, vanishing always occurs for any initial values. The next result
shows that spreading occurs for large values.
 \begin{thm} Suppose that $R_0^F(0)<1<R_0$. Then $h_\infty-g_\infty=\infty$ and $\liminf_{t\to
+\infty}||u(\cdot, t)||_{C([0, h(t)])}>0$ if $||u_{0}(x)||_{C([-h_0, h_0])}$ and $||v_{0}(x)||_{C([-h_0, h_0])}$ are sufficiently large.
\end{thm}
\bpf We  construct a vector $(\underline u , \underline v, \underline h)$ such that $u\geq \underline u$, $v\geq\underline v$ in $[-\underline h(t), \underline h(t)]\times [0, T_0]$, and also $g(t)\leq -\underline h(t)$, $h(t)\geq \underline h(t)$
in $[0, T_0]$. If we can choose $T_0$ such that $R_0^D((-\underline h(\sqrt{T_0}), \underline h(\sqrt{T_0}))>1$, then $h_\infty-g_\infty=\infty$.

We first consider the following eigenvalue problem
\begin{eqnarray}
\left\{
\begin{array}{lll}
-d\psi''-\frac 12\psi'=\mu_0 \psi,\; & 0<x<1,  \\
\psi (0)=\psi(1)=0. &
\end{array} \right.
\label{Bf11}
\end{eqnarray}
It is well known that the principal eigenvalue $\mu_0$ of this problem is simple;  the corresponding eigenfunction
$\psi(x)$ can be chosen to be  positive in $[0, 1)$ and $||\psi||_{L^\infty}=1$. It is also easy to see that $\mu_0>\frac 1{16d}$
 and $\psi'(x)<0$ in $(0, 1]$. Extending $\psi$ in $[0, 1]$ to an even function $\phi$ in $[-1, 1]$ yields
 \begin{eqnarray}
\left\{
\begin{array}{lll}
-d\phi''-\frac {{\textrm sgn}(x)}2\phi'=\mu_0 \phi,\; & -1<x<1,  \\
\phi(-1)=\phi(1)=0. &
\end{array} \right.
\label{Bh11}
\end{eqnarray}

We now construct a suitable lower solution to  \eqref{a3} and we define
$$\underline{h}(t)=\sqrt{t+\delta},\, t\geq 0, $$
$$\underline{u}(x,t)= \frac{M}{(t+\delta)^k}\phi(\frac{x}{\sqrt{t+\delta}}),\;  -\sqrt{t+\delta}\leq x\leq \sqrt{t+\delta}, \, t\geq 0,$$
$$\underline{v}(x,t)= \frac{M}{(t+\delta)^k}\phi(\frac{x}{\sqrt{t+\delta}}),\;  -\sqrt{t+\delta}\leq x\leq \sqrt{t+\delta}, \, t\geq 0,$$
where $\delta, M, T_0, k$ are chosen as follows :
$$0<\delta\leq \min\{1, h_0^2\},\ R_0^D((-\sqrt{T_0}, \sqrt{T_0}))>1,$$
$$k\geq \max\{\mu_0 +a_{11}(T_0+1),\, a_{22}(T_0+1)\},\ -2\mu M \phi'(1)>(T_0+1)^k.$$
Here we have used the assumption that $R_0>1$ and the fact that $R_0^D((-\sqrt{T_0}, \sqrt{T_0}))\to R_0$ as $T_0\to \infty$ by Lemma 3.4.

 Direct computations yield
\begin{eqnarray*}
& &\dfrac{\partial \underline u}{\partial t}-d \dfrac{\partial^2\underline u}{\partial x^2}+a_{11}\underline u-a_{12}\underline v\\
& &=-\frac{M}{(t+\delta)^{k+1}}[d \phi''+\frac {x}{2\sqrt{t+\delta}}\phi'+(k+(-a_{11}+a_{12})(t+\delta)) \phi]\\
& &\leq -\frac{M}{(t+\delta)^{k+1}}[d \phi''+\frac {sgn(x)}{2}\phi'+\mu_0\phi]\\
& &= 0,
\end{eqnarray*}
\begin{eqnarray*}
& &\dfrac{\partial \underline v}{\partial t}+a_{22} \underline v-G(\underline u(x,t)) \\
& &\leq -\frac{M}{(t+\delta)^{k+1}}[(k+(-a_{22}+G(\underline u)/\underline u)(t+\delta)) \phi]\\
& &\leq -\frac{M}{(t+\delta)^{k+1}}[k -a_{22}(T_0+1)]\\
& &\leq 0
\end{eqnarray*}
for all $0<t\leq T_0$ and $-\underline h<x<\underline h$.
$$\underline h'(t)+\mu \underline u_x(\sqrt{t+\delta},t)=\frac {1}{2\sqrt{t+\delta}}+\frac{\mu M}{(t+\delta)^{k+1/2}}\phi'(1)<0.$$
Then we have
\begin{eqnarray*}
\left\{
\begin{array}{lll}
\dfrac{\partial \underline u}{\partial t}\leq d \dfrac{\partial^2 \underline u}{\partial x^2}
-a_{11}\underline u +a_{12}\underline v,\; & -\underline h<x<\underline h,\, 0<t\leq T_0, \\
\dfrac{\partial \underline v}{\partial t}\leq
-a_{22}\underline v+G(\underline u(x,t)),\; &  -\underline h<x<\underline h,\, 0<t\leq T_0,\\
\underline {u}(x,t)=\underline{v}(x,t)=0,\; & x=\pm \underline h(t),\,0<t\leq T_0,\\
-\underline h_0=-\sqrt{\delta}\geq -h_0, & 0<t\leq T_0, \\
\underline h_0=\sqrt{\delta}\leq h_0, & 0<t\leq T_0\\
\underline h'(t)<-\mu \underline u_x(\sqrt{t+\delta},t),&0<t\leq T_0,\\
-\underline h'(t)>-\mu \underline u_x(-\sqrt{t+\delta},t),&0<t\leq T_0.
\end{array} \right.
\end{eqnarray*}
If $\underline u(x,0)=\frac{M}{\delta^{k}}\phi(\frac {x}{\sqrt{\delta}})<u_{0}(x)$ and
$\underline v(x,0)=\frac{M}{\delta^{k}}\phi(\frac {x}{\sqrt{\delta}})<v_{0}(x)$ in $[0, \sqrt{\delta}]$,
then using Lemma 2.5 yields  $h(t)\geq\underline h(t)$ and $g(t)\leq -\underline h(t)$
in $[0,T_0]$. In particular, $h(T_0)-g(T_0)\geq 2\underline h(T_0)=2\sqrt{T_0+\delta}\geq 2\sqrt{T_0}$. Noting that
\begin{eqnarray*}
R_0^F(T_0)&:=&R_0^D((g(T_0), h(T_0)))\geq R_0^D((-\underline h(T_0), \underline h(T_0)))\\
&\geq& R_0^D((-\sqrt{T_0}, \sqrt{T_0}))>1,
\end{eqnarray*}
 we then have $h_\infty-g_\infty=+\infty$ by Theorem 4.1.
 \epf

 \begin{thm} (Sharp threshold) Suppose that $R_0>1$,
with fixed $\mu$, $h_0$ and $(\phi,\psi)$ satisfying $(\ref{Ae1})$. Let $(u, v; g, h)$ be a solution
of $(\ref{a3})$ with $(u_0, v_0)=(\sigma \phi(x),\sigma \psi(x))$ for some $\sigma>0$. Then there exists $\sigma^*=\sigma^*(\phi,\psi)\in [0, \infty)$
 such that spreading occurs when $\sigma> \sigma^*$, and vanishing occurs when $0<\sigma\leq \sigma^*$.
\end{thm}
\bpf
It follows from Theorem 4.1 that spreading always occurs if $R_0^F(0)\geq 1$. Hence, in this
case we have $\sigma^*(\phi, \psi)=0$ for any $\phi$ and $\psi$.

For the remaining case $R_0^F(0)<1$,  define
$$\sigma^*:=\sup \{\sigma_0: h_\infty (\sigma\phi,\sigma \psi)<\infty \ \textrm{for}\ \sigma\in (0,\sigma_0]\}.$$
By Theorem 3.6, we
see that in this case vanishing occurs for all small $\sigma>0$; therefore, $\sigma^*\in
(0, \infty]$. On the other hand, it follows from Theorem 4.2 that in
this case spreading occurs for all large $\sigma$. Therefore, $\sigma^*\in
(0, \infty)$,  spreading occurs when $\sigma> \sigma^*$, and vanishing occurs when $0<\sigma< \sigma^*$ by Corollary 2.6.

We claim that vanishing occurs when $\sigma=\sigma^*$. Otherwise $h_\infty-g_\infty=\infty$ for
$\sigma=\sigma^*$. Since $R_0^F(t)\to R_0>1$ as $t\to \infty$,
 there exists $T_0>0$ such that $R_0^F(T_0)>1$. By
the continuous dependence of $(u, v; g, h)$ on its initial values, we can find
$\epsilon>0$ sufficiently small so that the solution of  (\ref{a3}) with $(u_0, v_0)=(\sigma^*-\epsilon) (\phi(x),\psi(x))$, denoted by
$(u_\epsilon, v_\epsilon; g_\epsilon, h_\epsilon)$ satisfies $R_0^F(T_0)>1$. This implies that spreading
occurs for $(u_\epsilon, v_\epsilon; g_\epsilon,h_\epsilon)$,
contradicting the definition of $\sigma^*$. This completes the proof.
 \epf

Similarly, if we consider $\mu$ instead of $u_0$ as a varying parameter, the following result holds;
 see also Theorem 4.4 in \cite{DL2}.
\begin{thm} (Sharp threshold) Suppose that $R_0>1$, with
fixed $h_0$, $u_0$ and $v_0$. Then there exists $\mu^*\in [0, \infty)$
 such that spreading occurs when $\mu> \mu^*$, and vanishing occurs when $0<\mu\leq \mu^*$.
\end{thm}

Next, we consider the asymptotic behavior of the solution to  \eqref{a3} when the spreading occurs.

\begin{thm} Suppose that $R_0>1$. If spreading occurs, then the solution of
 free boundary problem \eqref{a3} satisfies $\lim_{t\to +\infty} \ (u(x,t),v(x,t))=(u^*, v^*)$
uniformly in any bounded subset of $(-\infty, \infty)$, where $(u^*, v^*)$ is the unique positive equilibrium of
\eqref{aode1}.
\end{thm}
\bpf
(1) The limit superior of the solution

We recall that the comparison principle gives $(u(x,t),v(x,t))\leq (\overline u(t), \overline v(t))$
for $(x,t)\in [g(t), h(t)]\times (0,\infty)$, where
$(\overline u(t), \overline v(t))$ is the solution of the problem
    \begin{eqnarray}
\label{ode1}
\left\{
\begin{array}{lll}
&\overline u'(t)=-a_{11}\overline u(t)+a_{12}\overline v(t),&t>0, \\
&\overline v'(t)=-a_{22}\overline v(t)+G(\overline u(t)),&t>0,\\
&\overline u(0)=||u_0||_{L^\infty[-h_0, h_0]},\ \overline v(0)=||v_0||_{L^\infty[-h_0, h_0]}.&
\end{array} \right.
\end{eqnarray}
Since $R_0>1$, the unique positive equilibrium $(u^*, v^*)$ is globally stable for the ODE system (\ref{ode1}) and $\lim_{t\to\infty}(\overline
u(t), \overline v(t))= (u^*, v^*)$; therefore we  deduce
\begin{equation}\label{123}
\limsup_{t\to +\infty} \ (u(x,t), v(x,t))\leq (u^*, v^*)
\end{equation}
uniformly for $x\in (-\infty, \infty)$.

(2) The lower bound of the solution for a large time.

Note that $R_0>1$ and
$$\lim_{l\to \infty}\frac {G'(0)\frac {a_{12}}{a_{22}}}{a_{11}+d(\frac \pi {2l})^2}=R_0>1;$$
therefore, there is $L_0$ such that $\frac {G'(0)\frac {a_{12}}{a_{22}}}{a_{11}+d(\frac \pi {2L_0})^2}>1.$
This implies that the principal eigenvalue $\lambda_0^*$ of
\begin{eqnarray}
\left\{
\begin{array}{lll}
-d \psi_{xx}=-a_{11} \psi+{G'(0)\frac {a_{12}}{a_{22}}}\psi+\lambda_0^* \psi,\; &
x\in (-L_0, L_0),  \\
\psi(x)=0, &x=\pm L_0
\end{array} \right.
\label{B1k}
\end{eqnarray}
satisfies
$$\lambda_0^*=a_{11}+d(\frac \pi {2L_0})^2-G'(0)\frac {a_{12}}{a_{22}}<0.$$
Since  $h_\infty-g_\infty=\infty$ by assumption,  $h_\infty=g_\infty=\infty$ by Lemma 3.1. Thus,
 for any $L\geq L_0$, there exists $t_L>0$ such that $g(t)\leq -L$ and $h(t)\geq L$ for $t\geq t_L$.

Letting $\underline U=\delta \psi$ and $\underline V=G(\underline U)/a_{22}$,  we can choose $\delta$ sufficiently small such that $(\underline U, \underline V)$
satisfies
 \begin{eqnarray*}
\left\{
\begin{array}{lll}
\underline U_{t}\leq d \underline U_{xx}-a_{11}\underline U+a_{12}\underline V,\; &  -L_0<x<L_0,\ t>t_{L_0},  \\
\underline V_{t}=-a_{22}\underline V+G(\underline U),\; & -L_0<x<L_0,\ t>t_{L_0},  \\
\underline U(x,t)= 0,\ \quad & x=\pm L_0,\ t>t_{L_0},\\
\underline U(x,t_{L_0})\leq u(x,t_{L_0}),\ \underline V(x,t_{L_0})=v(x, t_{L_0}),
 & -L_0\leq x\leq L_0,
\end{array} \right.
\end{eqnarray*}
 meaning that $(\underline U, \underline V)$ is a lower solution of the solution $(u,v)$ in $[-L_0, L_0]\times [t_{L_0}, \infty)$.
We then have $(u, v)\geq (\delta \psi, G(\delta \psi)/a_{22})$ in $[-L_0, L_0]\times [t_{L_0}, \infty)$, which implies that the solution can not decay to zero.

(3) The limit inferior of the solution.

We extend $\psi(x)$ to $\psi_{L_0}(x)$ by defining $\psi_{L_0}(x):=\psi(x)$ for $-L_0\leq x\leq L_0$ and $\psi_{L_0}(x):=0$ for $x<-L_0$ or $x>L_0$.
Now for $L\geq L_0$, $(u,v)$ satisfies
\begin{eqnarray}
\left\{
\begin{array}{lll}
u_{t}=d u_{xx}-a_{11}u+a_{12}v,\; &  g(t)<x<h(t),\ t>t_L,  \\
v_{t}=-a_{22}v+G(u),\; & g(t)<x<h(t),\ t>t_L,  \\
u(x,t)=0, \quad & x=g(t)\, \textrm{or}\, x=h(t),\ t>t_L,\\
u(x,t_L)\geq \delta \psi_{L_0},\ v(x,t_L)\geq G(\delta \psi_{L_0})/a_{22},
 & -L\leq x\leq L;
\end{array} \right.
\label{fs1}
\end{eqnarray}
therefore, we  have $(u, v)\geq (\underline u, \underline v)$ in $[-L, L]\times [t_L, \infty)$,
where $(\underline u, \underline v)$ satisfies
\begin{eqnarray}
\left\{
\begin{array}{lll}
\underline u_{t}=d \underline u_{xx}-a_{11}\underline u+a_{12}\underline v,\; &  -L<x<L,\ t>t_L,  \\
\underline v_{t}=-a_{22}\underline v+G(\underline u),\; & -L<x<L,\ t>t_L,  \\
\underline u(x,t)= 0,\quad & x=\pm L,\ t>t_L,\\
\underline u(x,t_L)=\delta \psi_{L_0},\ \underline v(x,t_L)=G(\delta \psi_{L_0})/a_{22},
 & -L\leq x\leq L.
\end{array} \right.
\label{fs11}
\end{eqnarray}
The system (\ref{fs11}) is quasimonotone increasing; therefore,
it follows from the upper and lower solution method
 and the theory of monotone dynamical systems ( \cite{HS} Corollary 3.6) that
$\lim_{t\to +\infty} \ (\underline u(x,t), \underline v(x,t))\geq (\underline u_L(x), \underline v_L(x))$ uniformly in
$[-L, L]$, where $(\underline u_L, \underline v_L)$ satisfies
\begin{eqnarray} \label{fs12}\left\{
\begin{array}{lll}
-d \underline u''_{L}=-a_{11}\underline u_L+a_{12}\underline v_L,\; &  -L<x<L,  \\
-a_{22}\underline v_L+G(\underline u_L)=0,\; & -L<x<L,  \\
\underline u_L(x)=0, &x=\pm L
\end{array} \right.
\end{eqnarray}
and is the minimum upper solution $(\delta \psi_{L_0},\ G(\delta \psi_{L_0})/a_{22})$.

Now we give the monotonicity and show that if $0<L_1<L_2$, then $\underline u_{L_1}(x)\leq \underline u_{L_2}(x)$
 in $[-L_1, L_1]$. The result is derived by comparing the boundary
conditions and initial conditions in (\ref{fs11}) for $L=L_1$ and $L=L_2$.

Let $L\to \infty$.  By classical elliptic regularity theory and a diagonal procedure, it follows that  $(\underline u_{L}(x), \underline v_{L}(x))$ converges
uniformly on any compact subset of $(-\infty, \infty)$ to $(\underline u_\infty, \underline v_\infty)$ that is continuous on $(-\infty, \infty)$ and satisfies
\begin{eqnarray*} \left\{
\begin{array}{lll}
-d \underline u''_{\infty}=-a_{11}\underline u_\infty+a_{12}\underline v_\infty,\; &  -\infty<x<\infty,  \\
-a_{22}\underline v_\infty+G(\underline u_\infty)=0,\; & -\infty<x<\infty,  \\
\underline u_\infty(x)\geq \delta \psi_{L_0},\ \underline v_\infty(x)\geq G(\delta \psi_{L_0})/a_{22}, &-\infty<x<\infty.
\end{array} \right.
\end{eqnarray*}

Next, we observe that $\underline u_\infty (x)\equiv u^*$ and $\underline v_\infty (x)\equiv v^*$, which can be derived by considering the problem
$$-d w''=-a_{11}w+\frac{a_{12}}{a_{22}}G(w).$$
The uniqueness of the positive solution follows from the assumption on $G$ and the condition $R_0>1$.

Now for any given $[-M, M]$ with $M\geq L_0$, since that $(\underline u_{L}(x), \underline v_{L}(x))\to (u^*, v^*)$ uniformly in $[-M, M]$, which is the compact subset  of $(-\infty, \infty)$, as $L\to \infty$, we deduce that for any $\varepsilon >0$, there exists $L^*>L_0$ such that  $(\underline u_{L^*}(x), \underline v_{L^*}(x))\geq  (u^*-\varepsilon, v^*-\varepsilon)$ in $[-M, M]$. As above, there is $t_{L^*}$ such that $[g(t), h(t)]\supseteq [-L^*, L^*]$ for $t\geq t_{L^*}$.
Therefore, 
$$(u(x,t), v(x,t))\geq (\underline u(x,t), \underline v(x,t))\ \textrm{in}\ [-L^*, L^*]\times [t_{L^*}, \infty),$$ and 
$$\lim_{t\to +\infty} \ (\underline u(x,t), \underline v(x,t))\geq (\underline u_{L^*}(x), \underline v_{L^*}(x))\ \textrm{in}\ [-L^*, L^*].$$
Using the fact that $(\underline u_{L^*}(x), \underline v_{L^*}(x))\geq (u^*-\varepsilon, v^*-\varepsilon)$ in $[-M, M]$ gives 
 $$\liminf_{t\to +\infty} \ (u(x, t), v(x,t))\geq (u^*-\varepsilon, v^*-\varepsilon)\ \textrm{in}\ [-M, M].$$
 Since $\varepsilon>0$ is arbitrary, we then have $\liminf_{t\to +\infty} \ u(x, t)\geq u^*$ and $\liminf_{t\to +\infty} \ v(x,t)\geq v^*$ uniformly  in $[-M,M]$, which together with (\ref{123})
imply that $\lim_{t\to +\infty} \ u(x,t)=u^*$
and  $\lim_{t\to +\infty} \ v(x,t)=v^*$ uniformly in any bounded subset of $(-\infty, \infty)$.
\epf

Combining Remark 4.1, Theorem 4.2 and Theorem 4.5, we immediately obtain the following
spreading-vanishing dichotomy:
\begin{thm} Suppose that $R_0>1$.
Let $(u(x, t), v(x, t); g(t), h(t))$ be the solution of  free boundary problem \eqref{a3}.
Then, the following alternatives hold:

Either
\begin{itemize}
\item[$(i)$] {\rm Spreading:} $h_\infty-g_\infty =+\infty$ and $\lim_{t\to +\infty} \ (u(x, t), v(x, t))=(u^*, v^*)$
uniformly in any bounded subset of $(-\infty, \infty)$; \end{itemize}

or
\begin{itemize}
\item[$(ii)$] {\rm Vanishing:} $h_\infty -g_\infty \leq h^*$ with $\frac {G'(0)\frac {a_{12}}{a_{22}}}
{a_{11}+d(\frac \pi {h^*})^2}=1$ and $\lim_{t\to +\infty} \ (||u(\cdot, t)||_{C([g(t),h(t)])}+||v(
\cdot, t)||_{C([g(t),h(t)])})=0$.
\end{itemize}
\end{thm}

\section{Discussion}

In this paper, a free boundary problem is used to describe the expanding of  bacteria in a man-environment-man
epidemic model in a one-dimensional habitat. We take into account the spreading and vanishing of the bacteria.
Here, vanishing implies not only that the bacteria disappear eventually, but also that the infected habitat is limited, and
spreading means the existence of the bacteria in the long run with an uncontrollable infected environment.
Sufficient conditions for the bacteria spreading or vanishing are given.

Compared with existing work described by reaction-diffusion systems \eqref{Aa1} in \cite{CM}
or established by travelling waves and entire solutions in \cite{TZ, WW, WW2, WS, XZ, ZW} ,
our model \eqref{a3} provides a different way to understand the expanding process of bacteria.
It is well-known that for the ODE system \eqref{aode1}, the basic reproduction number
$R_0(:=\frac {G'(0)a_{12}}{a_{11}a_{22}})$ determines
whether the bacteria die out ($R_0<1$) or remain endemic ($R_0>1$). However, in our problem \eqref{a3},
the infected habitat is changing with time $t$; therefore, we  introduced the basic reproduction number
$$R^F_0(t):=\frac {G'(0)\frac {a_{12}}{a_{22}}}{a_{11}+d(\frac \pi {h(t)-g(t)})^2},$$
which depends on the habitat $(g(t), h(t))$, the diffusion rate $d$ and the coefficients in \eqref{a3}.
We showed that $R^F_0(t)\leq R_0$ and $R^F_0(t)\to R_0$ if $(g(t), h(t))\to (-\infty, +\infty)$ as $t\to \infty$.
Furthermore, if $R_0\leq 1$, the bacteria are always vanishing (Theorem 3.5). The result is the same as that for the corresponding
ODE system \eqref{aode1}.  However, if $R^F_0(t_0)\geq 1$ for some $t_0\geq 0$, the bacteria
are spreading (Theorem 4.1 and Remark 4.1). For the case $R^F_0(0)<1<R_0$, the spreading or vanishing of the bacteria depends
on the initial size of  bacteria (Theorem 4.3 ), or the ratio ($\mu$) of the expansion speed of the free boundary
and the population gradient at the expanding fronts (Theorem 4.4).

Ecologically, our main results reveal that if the multiplicative factor of the infectious bacteria is small, the bacteria will
die out eventually and the humans are safe. Otherwise, the spreading or vanishing of the bacteria depends on
the initial infected habitat, the diffusion rate, and other factors. In particular, the initial number of  bacteria plays a
 key role. A large initial number can induce the spreading of bacteria easily.
 A similar result obtained for an invasive species has been supported by substantial  empirical evidence; see \cite{DL}.
    Therefore, we hope our model and theoretical results can be used to provide better prediction and prevention of
  infecting bacteria.

\end{document}